\documentclass{amsart}
\usepackage[latin1]{inputenc}
\usepackage{amssymb}
\usepackage[usenames,dvipsnames,dvips]{color}

\newtheorem{thm}{Theorem}
\newtheorem{cor}[thm]{Corollary}
\newtheorem{lem}[thm]{Lemma}
\newtheorem{prop}[thm]{Proposition}

\newtheorem{conj}[thm]{Conjecture}
\theoremstyle{definition}
\newtheorem{defn}[thm]{Definition}
\theoremstyle{remark}
\newtheorem{rem}[thm]{Remark}
\numberwithin{equation}{section}

\newcommand{\To}{\longrightarrow}

\begin{document}

\title{Compact spaces that do not map onto finite products}

\author{Antonio Avil\'{e}s}

\subjclass[2000]{}

\address{A. Avil\'{e}s, University of Paris 7, UMR 7056, 2, Place Jussieu, Case 7012, 75251 Paris (France)} \email{aviles@logique.jussieu.fr, avileslo@um.es}

\thanks{The author was supported by a Marie Curie
Intra-European Felloship MCEIF-CT2006-038768 (E.U.) and research
project MTM2005-08379 (MEC and FEDER)}

\begin{abstract}
We provide examples of nonseparable compact spaces with the
property that any continuous image which is homeomorphic to a
finite product of spaces has a maximal prescribed number of
nonseparable factors.

\end{abstract}

\maketitle

\section{Introduction}

The motivation of this work are several problems coming from
\cite{AvilesCountableProducts}, \cite{fiberorders} and
\cite{Mardesic}, dealing with the possibility of mapping a finite
product of spaces onto a another product space with more factors.

Consider $B$ to be the closed euclidean ball of a nonseparable
Hilbert space, endowed with its weak topology. It was proven in
\cite{fiberorders} that $B^2$ is not homeomorphic to $B$, and the
natural question arises whether $B^2$ is at least a continuous
image of $B$. We shall prove in this note that the answer to this
question is also negative. More generally,

\begin{thm}\label{Ballindecomposable}
Let $n<m$ be natural numbers and suppose $f:B^n\To
X_1\times\cdots\times X_m$ is an onto continuous map. Then there
exists $i\leq m$ such that $X_i$ is metrizable.
\end{thm}

We shall also provide an alternative proof of the fact shown in
\cite{fiberorders} that if $B^n$ is homeomorphic to $L^m$ for some
$L$ and some $m$, then $m$ divides $n$. These kind of properties
are also proven in \cite{fiberorders} for some spaces of
probability measures, but we shall show that the methods in this
paper do not apply to these spaces.

The second problem deals with linearly ordered spaces. The
following result was first obtained by Treybig~\cite{Treybig},
though there exists a shorter proof by Bula, Debski and
Kulpa~\cite{BulDebKul}:

\begin{thm}[Treybig]
Let $L$ be a linearly ordered compact space, and $X_0$ and $X_1$
two infinite compact spaces. If there is a continuous surjection
$f:L\To X_0\times X_1$, then both $X_0$ and $X_1$ are metrizable.
\end{thm}

Marde\v{s}i\'{c}~\cite{Mardesic} tried to generalize this theorem
to higher dimensions, proposing the following conjecture:

\begin{conj}[Marde\v{s}i\'{c}]
Let $L_1,\ldots, L_n$ be linearly ordered compact spaces, let
$X_0,\ldots,X_n$ be infinite compact spaces and let $f:\prod_1^n
L_i\To \prod_0^n X_j$ be an onto continuous map. Then there exist
$i,j\leq n$, $i\neq j$ such that $X_i$ and $X_j$ are metrizable.
\end{conj}

He proved \cite{Mardesic} that the conjecture holds under the
assumption that all spaces $X_i$ are separable. Our methods enable
us to obtain the following:

\begin{thm}\label{orderindecomposable}
Let $L_1,\ldots, L_n$ be linearly ordered compact spaces,let
$X_0,\ldots,X_n$ be infinite compact spaces and let $f:\prod_1^n
L_i\To \prod_0^n X_j$ be an onto continuous map. Then there exist
$0\leq i,j\leq n$, $i\neq j$ such that $X_i$ and $X_j$ are
separable.
\end{thm}

Notice that Marde\v{s}i\'{c}'s partial answer to the conjecture
and our own seem completely unrelated since his hypothesis is
stronger than our conclusion. However, in the case when $n=2$ both
results can be combined to conclude that at least one factor must
be always metrizable:

\begin{cor}
Let $L_1$ and  $L_2$ be linearly ordered compact spaces,let $X_0$,
$X_1$ and $X_2$ be infinite compact spaces and let $f: L_1\times
L_2\To X_0\times X_1\times X_2$ be an onto continuous map. Then
there exists $i\in\{0,1,2\}$ such that $X_i$ is metrizable.
\end{cor}

Proof: By Theorem~\ref{orderindecomposable}, there are two factors
which are separable, say that $X_0$ and $X_1$ are separable. Let
$Y_3$ be an infinite quotient space of $X_3$ of countable weight.
Then there is a quotient $L_1\times L_2\To X_0\times X_1\times
Y_3$ where all factors are separable, so by Marde\v{s}i\'{c}'s
result~\cite{Mardesic} at least two of them are metrizable, so
either $X_0$ or $X_1$ is metrizable.$\qed$

The third problem deals with the spaces $\sigma_n(\Gamma) = \{x\in
2^\Gamma : |supp(x)|\leq n\}$. In our previous work
\cite{AvilesCountableProducts} we studied the homeomorphic
classification of finite and countable products of these spaces.
In this paper, we shall determine when such a product is a
continuous image of another one.

We want to express our gratitude to Stevo Todorcevic for calling
our attention to the work of Marde\v{s}i\'{c}.

\section{Indecomposability properties (I)}

All along this paper, we work with compact Hausdorff topological
spaces. When we talk about compact spaces, the Hausdorff $T_2$
separation axiom is implicitly assumed.

\begin{defn}
Let $X$ be a compact space. A pseudoclopen of $X$ is a pair
$a=(a[0],a[1])$ such that $a[0]$ and $a[1]$ are open subsets of
$X$ and $\overline{a[0]}\subset a[1]$.
\end{defn}

Notice that every clopen set $c$ can be identified with a
pseudoclopen $(c,c)$. Conversely, if $K$ is a totally disconnected
compact space, then every pseudoclopen $a$ is interpolated by a
clopen set $c$, meaning $a[0]\subset c \subset a[1]$. The notion
of pseudoclopen substitutes the notion of clopen sets in general
(not totally disconnected) compact spaces.

\begin{defn}
An uncountable family $\mathcal{F}$ of sets will be called a
Knaster-disjoint family if every uncountable subfamily
$\mathcal{G}\subset\mathcal{F}$ contains two disjoint elements. An
uncountable family $\mathcal{F}$ of pseudoclopens of $X$ will be
called Knaster-disjoint if $\{a[1] : a\in\mathcal{F}\}$ is
Knaster-disjoint family of sets.
\end{defn}

The reason for the name Knaster-disjoint is because of the well
known Knaster's chain condition: \emph{Every uncountable family of
nonempty opens sets contains an uncountable family in which every
two elements have nonempty intersection}. Thus, for a completely
regular space, the failure of Knaster's condition is equivalent to
the existence of an uncountable Knaster-disjoint family of
strongly nonempty pseudoclopens (we call a pseudoclopen $a$
strongly nonempty if $a[0]\neq\emptyset$).

\begin{defn}
Let $X$ be a compact space, and $n$ a natural number. We say that
$X$ has property $I_n$ if for every $n+1$ many Knaster-disjoint
families $\mathcal{F}_0,\ldots,\mathcal{F}_n$ of pseudoclopens,
there exist uncountable subfamilies $\mathcal{G}_i\subset
\mathcal{F}_i$ for $i=0,\ldots,n$ such that for every choice
$a_i\in \mathcal{G}_i$ for $i=0,\ldots,n$ we have
$a_0[0]\cap\cdots\cap a_n[0]=\emptyset$.
\end{defn}

\begin{prop}\label{imageIn}
Let $X$ be a compact space with property $I_n$, and $f:X\To Y$ a
continuous onto map. Then $Y$ has property $I_n$.
\end{prop}

Proof: Let $\mathcal{F}_0,\ldots,\mathcal{F}_n$ be families of
Knaster-disjoint pseudoclopens in $Y$. For every $i\leq n$, let
$\mathcal{F}'_i = \{(f^{-1}(a[0]), f^{-1}(a[1])) : a\in
\mathcal{F}_i\}$. Then $\mathcal{F}'_i$ is a Knaster-disjoint
family of pseudoclopens in $X$. Because $X$ has property $I_n$, we
can find uncountable subfamilies $\mathcal{G}_i\subset
\mathcal{F}_i$ such that whenever $a_i\in G_i$, then
$f^{-1}(a_0[0])\cap\cdots\cap f^{-1}(a_n[0])=\emptyset$, which
implies that $a_0[0]\cap\cdots\cap a_n[0]=\emptyset$ because $f$
is onto. $\qed$

\begin{prop}\label{onecccIn}
Let $X_0,\ldots,X_n$ be compact spaces such that
$X=X_0\times\cdots\times X_n$ has property $I_n$. Then there
exists $i\leq n$ such that $X_i$ satisfies Knaster's condition.
\end{prop}

Proof: For a set $s\subset X_i$, we write $u_i(s)
=\{(x_0,\ldots,x_n)\in X : x_i\in s\}$. Also, for a pseudoclopen
$a$ of $X_i$, $u_i(a) = (u_i(a[0]),u_i(a[1]))$ is a pseudoclopen
of $X$. Now, suppose no $X_i$ satisfies Knaster's condition. Then,
for every $i$ we can find an uncountable Knaster-disjoint family
$\mathcal{F}_i$ of strongly nonempty pseudoclopens in $X_i$.
Associated with them, we have  Knaster-disjoint families
$\mathcal{F}'_i = \{u_i(a) : a\in \mathcal{F}_i\}$ of
pseudoclopens in $X$. These families contradict property $I_n$,
because $a_0[0]\cap\cdots\cap a_n[0]\neq\emptyset$ whenever
$a_i\in\mathcal{F}'_i$.$\qed$

We introduce now some operations and relations. For two
pseudoclopens $a$ and $b$, we write:
\begin{itemize}
\item $b\subset a$ if $b[0]\subset a[0]$ and $b[1]\subset a[1]$,
\item $b\prec a$ if $a[0]\subset b[0]\subset b[1]\subset a[1]$ (we
say that $b$ is finer than $a$). \item $a\cup b = (a[0]\cup
b[0],a[1]\cup b[1])$ (notice that this a new pseudoclopen).
\end{itemize}

\begin{defn}A family $\mathcal{B}$ of pseudoclopens of $X$ will be called
a basis if for every pseudoclopen $a$ of $X$ there exists
$b_1,\ldots,b_m$ finitely many pseudoclopens from $\mathcal{B}$
such that $b_1\cup\cdots\cup b_m\prec a$.
\end{defn}

\begin{lem}
If $\mathcal{B}$ is basis for the topology of a compact space $X$,
then the family of all pseudoclopens $a$ such that
$a[0]\in\mathcal{B}$ and $a[1]\in\mathcal{B}$ constitutes a basis
of pseudoclopens.
\end{lem}

\begin{lem}\label{clopenbasis}
If $X$ is a totally disconnected compact space and $\mathcal{B}$
is a basis for the topology of $X$ consisting of clopen sets, then
$\mathcal{B'} = \{(c,c) : c\in\mathcal{B}\}$ is a basis of
pseudoclopens in $X$.
\end{lem}

\begin{lem}\label{basislemma}
In the definition of property $I_n$ we may assume that the
families $\mathcal F_i$ are subfamilies of some basis of
pseudoclopens for the space.
\end{lem}

Proof: Let $\mathcal{F}_0,\ldots,\mathcal{F}_n$ be arbitrary
Knaster-disjoint families of pseudoclopens of $X$, and let
$\mathcal{B}$ be a basis of pseudoclopens. Because $\mathcal{B}$
is a basis, we can suppose (by passing to finer pseudoclopens)
that every $a\in\mathcal{F}_i$ is a finite union of elements of
$\mathcal{B}$. Because we are looking for uncountable subfamilies,
we can suppose without loss of generality that there exist natural
numbers $m_0,\ldots,m_n$, such that every element $a\in
\mathcal{F}_i$ is a union of exactly $m_i$ elements of
$\mathcal{B}$, and we write it in the form $a =
a_i^1\cup\cdots\cup a^{m_i}_i$ where each $a^j_i$ is an element of
$\mathcal{B}$. For every $i\leq n$ and every $j\leq m_i$, let
$\mathcal{F}_i^j = \{a_i^j : a\in \mathcal{F}_i\}$. For every
choice of numbers $j=(j_0\ldots,j_n)$ with $j_i\leq m_i$, we can
apply our hypothesis to the families
$\mathcal{F}_0^{j_0},\ldots,\mathcal{F}_n^{j_n}$, and this implies
that are uncountable subfamilies
$\mathcal{G}_i^{j}\subset\mathcal{F}_i$ such that whenever $a\in
\mathcal{G}_i^j$, we have that $a_0^{j_0}[0]\cap\cdots\cap
a_n^{j_n}[0]=\emptyset$. There are only finitely many choices for
the tuple $j=(j_0,\ldots,j_n)$, say that we call them
$j^{(1)},\ldots, j^{(k)}$, so we can find the subfamilies
$\mathcal{G}_i^j$ successively one after another and satisfying
$\mathcal{G}_i^{j^{(r)}}\supset \mathcal{G}_i^{j^{(r+1)}}$. The
uncountable families $\mathcal{G}_i^{j^{(k)}}\subset
\mathcal{F}_i$ satisfy the desired conclusion.$\qed$

The following fact is well known: it is the key property behind
the fact that, unlike the countable chain condition, Knaster's
condition is productive~\cite{Szpilrajn}:

\begin{lem}\label{Knasterrectangles}
Let $\{a^0\times a^1 : a\in\mathcal{F}\}$ be a Knaster-disjoint
family of subsets of $X\times Y$ consisting of rectangles. Then
one of the two families $\{a^0: a\in\mathcal{F}\}$ or $\{a^1 :
a\in\mathcal{F}\}$ contains an uncountable Knaster-disjoint
subfamily.
\end{lem}

\begin{prop}\label{productivity}
Let $X$ and $Y$ be compact spaces with property $I_n$ and $I_m$
respectively. Then $X\times Y$ has property $I_{n+m}$.
\end{prop}

Proof: Let $\mathcal{F}_0,\ldots,\mathcal{F}_{n+m}$ be uncountable
Knaster-disjoint families of psuedoclopens of $X\times Y$. By
Lemma~\ref{basislemma}, we can suppose that every
$a\in\mathcal{F}_i$ is of the form $a = a^0\times a^1 =
(a^0[0]\times a^1[0],a^0[1]\times a^1[1])$ with $a^0$ and $a^1$
pseudoclopens of $X$ and $Y$. By Lemma~\ref{Knasterrectangles},we
can also suppose that for every $i\leq n+m$ there exists
$j(i)\in\{0,1\}$ such that the family $\{a^{j(i)} :
a\in\mathcal{F}_i\}$ is Knaster-disjoint. By elementary
cardinality reasons, either there exists $S\subset
\{0,\ldots,n+m\}$ with $|S|=n+1$ such that $(\forall i\in S)
(j(i)=0)$, or else there exists $T\subset \{0,\ldots,n+m\}$ with
$|T|=m+1$ such that $(\forall i\in T) (j(i)=1)$. In the first
case, we finish the proof appealing to property $I_n$ of $X$, and
in the second case appealing to property $I_m$ of $Y$.$\qed$

\begin{prop}
Let $X$ be a compact space and $d$, $n$ and $q$ natural numbers.
Suppose that $X^d$ has property $I_n$ and $(q+1)d>n$. Then $X$ has
property $I_q$
\end{prop}

Proof: Let $\mathcal{F}_0,\ldots,\mathcal{F}_q$ be
Knaster-disjoint families of pseudoclopens of $X$. For every
$i\in\{1,\ldots,d\}$ let $p_i:X^d\To X$ be the projection on the
$i$-th coordinate. Let
$$\mathcal{F}_j^i = \{ (p_i^{-1}a[0],p_i^{-1}a[1]) :
a\in\mathcal{F}_j\}.$$ These are $(q+1) d$ many Knaster-disjoint
families of pseudoclopens of $X^d$. Since $(q+1)d>n$, and $X^d$
has property $I_n$ it follows that there are uncountable
subfamilies $\mathcal{G}_j^i\subset \mathcal{F}_j^i$ such that
$$\bigcap_{i=1}^d\bigcap_{j=0}^q c_j^i[0]=\emptyset, \text{ whenever } c_j^i\in\mathcal{G}_j^i.$$
We claim that there exists $i\in\{1,\ldots,d\}$ such that
$b_0[0]\cap\cdots\cap b_q[0]=\emptyset$ whenever
$b_0\in\mathcal{G}^i_0,\ldots,b_q\in\mathcal{G}^i_q$. This claim
concludes the proof because then the families $\mathcal{G}_j =
\{a\in\mathcal{F}_j : (p_i^{-1}a[0],p_i^{-1}a[1])\in
\mathcal{G}_j^i\}$ are the uncountable subfamilies of the families
$\mathcal{F}_j$ that we look for. The claim is proved by
contradiction. It it was false we could find for every $i$
elements $b^i_0\in\mathcal{G}^i_0,\ldots,b^i_q\in\mathcal{G}^i_q$
such that $b^i_0[0]\cap\cdots\cap b^i_q[0]\neq\emptyset$. Since a
clopen from a family $\mathcal{G}^i_j$ depends only on the
coordinate $i$, this implies that $\bigcap_{i=1}^d\bigcap_{j=0}^q
b^i_j[0]\neq\emptyset$, which is a contradiction$\qed$

\begin{cor}\label{roots}
Let $K$ be a compact space with property $I_n$ but not $I_{n-1}$.
If $K$ is homeomorphic to $X^d$ for some space $X$, then $d$
divides $n$.
\end{cor}

Proof: If $d$ does not divide $n$, then there exists an integer
$q$ such that $q<\frac{n}{d}<q+1$. The previous proposition yields
that $X$ has property $I_q$, so by Proposition~\ref{productivity}
$X^d\approx K$ has $I_{qd}$, and therefore $I_{n-1}$ because
$qd<n$, so $qd\leq n-1$.$\qed$

\section{The euclidean ball}

Let $\Gamma$ be an uncountable set, let $$B=B(\Gamma) =
\left\{x\in [-1,1]^\Gamma : \sum_{\gamma\in\Gamma}|x_\gamma|\leq
1\right\}.$$

We consider this as a compact space endowed with the pointwise
topology. This $B(\Gamma)$ is actually homeomorphic to the ball of
the Banach space $\ell_p(\Gamma)$ in the weak topology for
$1<p<\infty$ and also to the dual ball of $c_0(\Gamma)$ in the
weak$^\ast$ topology.

\begin{thm}\label{Ballirreducible}
The space $B$ has property $I_1$.
\end{thm}

We observe that this result implies
Theorem~\ref{Ballindecomposable}. Let $n<m$ and $f:B^n\To
X_1\times\cdots\times X_n$ a continuous surjection. By
Proposition~\ref{imageIn} and Proposition~\ref{productivity},
$X_1\times\cdots\times X_n$ has property $I_n$ and by
Proposition~\ref{onecccIn} there exists $i\leq m$ such that $X_i$
has Knaster's condition, therefore also the countable chain
condition (every disjoint family of opens sets is countable). But
$B$ is an Eberlein compact, so its continuous image $X_i$ is also
Eberlein compact \cite{BenRudWag} and a result of
Rosenthal~\cite{Rosenthal} establishes that the countable chain
condition implies countable weight for an Eberlein compact.

It is easy to see that $B_+(\Gamma) = B(\Gamma)\cap [0,1]^\Gamma$
maps continuously onto $B(\Gamma)$; an onto continuous map
$B_+(\Gamma\times 2)\To B(\Gamma)$ is given by $x\mapsto
(x_{(\gamma,0)}-x_{(\gamma,1)})_{\gamma\in \Gamma}$. For
notational simplicity, we write $\lambda_n = \frac{1}{2^{n+1}}$.
The following totally disconnected compact space maps onto
$B_+(\Gamma)$:

$$L = \left\{x\in \{0,1\}^{\Gamma\times\omega} :
\sum_{(\gamma,n)\in\Gamma\times\omega}\lambda_n x_{\gamma,n}\leq
1\right\}.$$

The surjection $g:L\To B_+(\Gamma)$ is given by $g(x)_\gamma =
\sum_{n<\omega}\lambda_n x_{\gamma,n}$. By
Proposition~\ref{imageIn}, Theorem~\ref{Ballirreducible} follows
from the following:

\begin{thm}\label{LhasI1}
The space $L$ has property $I_1$.
\end{thm}

Proof: We consider $\mathcal{B}$ the basis of clopen subsets of
$L$ consisting of the sets of the form
$$a_U^V = \{x\in L : \forall (\gamma,n)\in U\ x_{\gamma,n}=1,\ \forall
(\gamma,n)\in V\ x_{\gamma,n}=0\},$$

where $U$ and $V$ are two disjoint finite subsets of
$\Gamma\times\omega$. It will convenient to use the following
notations: given $a=a_U^V\in\mathcal{B}$, we shall call $U(a)=U$
and $V(a)=V$. Also, for a finite set $U\subset
\Gamma\times\omega$, we call $$\sigma(U) = \sum_{(\gamma,n)\in
U}\lambda_n = \sum_{n<\omega}\lambda_n|U\cap \Gamma\times\{n\}|.$$
Notice the fundamental property that if $U\cap U'=\emptyset$ then
$\sigma(U\cup U')=\sigma(U)+\sigma(U')$. We need to know when two
elements of $\mathcal{B}$ are disjoint:

(DC) Let $a,b\in\mathcal{B}$. Then $a\cap b=\emptyset$ if and only
if one of the three following conditions holds: either $U(a)\cap
V(b)\neq \emptyset$ or $V(a)\cap U(b)\neq\emptyset$ or
$\sigma(U(a)\cup U(b))>1$.

Let $\mathcal{F}_0$ and $\mathcal{F}_1$ be two Knaster-disjoint
families of clopen sets from $\mathcal{B}$. By
Lemma~\ref{clopenbasis} and Lemma~\ref{basislemma}, it is enough
that we check property $I_1$ on such two families of clopen sets.

By the familiar $\Delta$-system lemma (cf. \cite[Theorem
9.18]{Jechsbook}), we can assume that the families $\{U(a) :
a\in\mathcal{F}_0\}$, $\{U(a) : a\in\mathcal{F}_1\}$, $\{V(a) :
a\in\mathcal{F}_0\}$ and $\{V(a) : a\in\mathcal{F}_1\}$ are
$\Delta$-systems of roots $R_0$, $R_1$, $S_0$ and $S_1$
respectively. We write $U(a) = R_i\cup U_i(a)$ and $V(a) = S_i\cup
V_i(a)$ for $a\in\mathcal{F}_i$, separating the root and the
disjoint part of the $\Delta$-system in such a way that $\{U_i(a)
: a\in\mathcal{F}_i\}$ and $\{V_i(a) : a\in\mathcal{F}_i\}$ are
disjoint families of finite sets for $i=0,1$. By a further
refinement we can also suppose that the whole family
$$(\star)\ \{U_0(a),U_1(b),V_0(a),V_1(b)
: a\in\mathcal{F}_0, b\in\mathcal{F}_1\}\cup \{R_0\cup R_1\cup
S_0\cup S_1\}$$ is a disjoint family.

The number $\sigma(U)$ is always a rational number, so we can also
suppose that $\sigma(U_0(a))= q_0$ and $\sigma(U_1(b))= q_1$ are
rational numbers independent of $a\in\mathcal{F}_0$ and
$b\in\mathcal{F}_1$.

\emph{Claim}: $\sigma(R_i)+ 2q_i>1$ for $i=0,1$.

Proof of the claim: Since $\mathcal{F}_i$ is a Knaster-disjoint
family, we can pick $a,b\in\mathcal{F}_i$ two different elements
such that $a\cap b=\emptyset$. Thus, one of the three alternatives
of the disjointness criterion (DC) must hold. But the two first
alternatives are impossible. For example, the disjointness of the
family $(\star)$ above implies that $U(a)\cap V(b) = R_i \cap S_i
\subset U(a)\cap V(a)=\emptyset$. Therefore, the third alternative
holds:
$$1<\sigma(U(a)\cup U(b)) = \sigma(R_i) +\sigma(U_i(a))
+\sigma(U_i(b))= \sigma(R_i) + 2q_i.$$ We finish the proof by
showing that, after all these refinements, $a\cap b=\emptyset$
whenever $a\in\mathcal{F}_0$ and $b\in\mathcal{F}_1$. Using again
the disjointness criterion (DC) we prove that $\sigma(U(a)\cup
U(b))>1$.
$$\sigma(U(a)\cup U(b)) = \sigma(R_0\cup R_1) +
\sigma(U_0(a))+\sigma(U_1(b)) = \sigma(R_0\cup R_1) + q_0 + q_1.$$
Say that $q_i = \min(q_0,q_1)$, then by Claim (A),
$$\sigma(U(a)\cup U(b)) = \sigma(R_0\cup R_1) + q_0 + q_1 \geq \sigma(R_i) + 2q_i > 1.\qed$$

\bigskip

\begin{cor}[Avil\'{e}s, Kalenda]\label{rootsB}
Let $X$ be a compact space and $m,n$ natural numbers such that
$B^n$ is homeomoprhic to $X^m$. Then $m$ divides $n$.
\end{cor}

Proof: Apply the preceding theorem and
Corollary~\ref{roots}.$\qed$

\section{Remarks about spaces $P(K)$}

Given a compact space $K$, we denote by $P(K)$ the space of Radon
probability measures on $K$ endowed with the weak$^\ast$ topology.
Results analogous to Corollary~\ref{rootsB} are proven in
\cite{fiberorders} for certain spaces of probability measures
(like $P(\sigma_n(\Gamma))$ and $P(\sigma_1(\Gamma)^n)$), so it is
a natural question whether such spaces have property $I_1$. We
show in this section that property $I_n$ on $P(K)$ imposes very
restrictive conditions on $K$.

\begin{prop}
Let $X$ be a compact space which contains $n$ many open subsets
whose closures are pairwise disjoint and fail the countable chain
condition. Then $P(X)$ maps continuously onto $B^n$, and in
particular $P(X)$ does not have property $I_{n-1}$.
\end{prop}

Proof: Let $V_1,\ldots, V_n$ be open subsets of $X$ with
$\overline{V}_i\cap\overline{V}_j=\emptyset$ for $i\neq j$, and
for every $i$, let $\mathcal U_i$ be an uncountable disjoint
family of nonempty open subsets of $V_i$. For every
$u\in\bigcup_{i=1}^n \mathcal U_i$ let $h_u:X\To [0,1]$ be a
continuous function such that $h_u(X\setminus u)=0$ and
$\max\{h_u(x) : x\in u\}=1$. For every $i\leq n$ let also
$g_i:X\To [0,1]$ be a continuous function such that $g_i(x)=0$ if
$x\in V_i$ and $g_i(x)=1$ if $x\in V_j$ for some $j\neq i$.
Let $$\chi_n(t)= \left\{%
\begin{array}{ll}
    n & \text{ if } t\geq 1-\frac{1}{n}, \\
    \frac{1}{1-t} & \text { if } t<1-\frac{1}{n}.
\end{array}
\right.$$ For $\mu$ a Radon measure on $X$ and $\phi:X\To
\mathbb{R}$ a continuous function on $X$, we put $\mu(\phi)= \int
\phi d\mu$. We define $f:P(X)\To [0,1]^{\mathcal U_1}
\times\cdots\times [0,1]^{\mathcal U_n}$ as follows: $$f(\mu)_u =
\chi_n(\mu(g_i))\cdot\mu(h_u),\ \ u\in\mathcal U_i.$$ For a set
$\Gamma$, remember that $B_+(\Gamma) = \{x\in[0,1]^\Gamma :
\sum_{\gamma\in\Gamma} x_\gamma\leq 1\}$. $B(\Gamma)$ is a
continuous image of $B_+(\Gamma)$, and we shall show that $f(P(X))
= B_+(\mathcal U_1)\times\cdots\times B_+(\mathcal U_n)$.

 First, we check that $f(P(X)) \subset
B_+(\mathcal U_1)\times\cdots\times B_+(\mathcal U_n)$. For fixed
$\mu\in P(X)$ and $i\leq n$, $$\sum_{u\in\mathcal U_i} f(\mu)_u =
\chi_n(\mu(g_i))\sum_{u\in\mathcal U_i}\mu(h_u)\leq
\frac{1}{1-\mu(g_i)}\sum_{u\in\mathcal U_i}\mu(h_u)\leq 1$$
because $g_i$ has disjoint support from all $h_u$'s, so
$\mu(g_i)+\sum_{u\in\mathcal U_i}\mu(h_u)\leq 1$.

Conversely, we prove now that any $x\in B_+(\mathcal
U_1)\times\cdots\times B_+(\mathcal U_n)$ belongs to $f(P(X))$.
For every $i\leq n$ let $\xi_i\in\overline{V}_i$,
$\xi_i\not\in\bigcup\{u : u\in \mathcal{U}_i\}$, and for every
$u\in\mathcal U_i$ let $\zeta_u\in u$ such that $h_u(\zeta_u)=1$.
We define $\mu\in P(X)$ a discrete probability measure on $X$ such
that $\mu\{\zeta_u\} = \frac{x_u}{n}$ for every $u$, $\mu\{\xi_n\}
= \frac{1-\sum_{u\in\mathcal U_n}x_u}{n}$. We have that $\mu(h_u)
= \frac{x_u}{n}$ for $u\in\mathcal{U}_n$, $\mu(g_i) =
1-\frac{1}{n}$, and $f(\mu) = x$.$\qed$

Remember that $\sigma_1(\Gamma)$ is the one point compactification
of a discrete set $\Gamma$. $P(\sigma_1(\Gamma))$ is homeomorphic
to $B_+(\Gamma)$, a continuous image of $B(\Gamma)$, so it has
property $I_1$ by Theorem~\ref{Ballirreducible}. Also, if $K_n$ is
a discrete union of $n$ many disjoint copies of
$\sigma_1(\Gamma)$, then $P(K_n)$ has property $I_n$, because it
is a continuous image of $P(\sigma_1(\Gamma))^n\times
B_+(\{1,\ldots,n\})$. We may ask if there is some sufficent
conditions on $L$ so that $P(L)$ has property $I_n$.

\section{Indecomposability properties (II)}

A family $\mathcal{F}$ of pseudoclopens will be called disjoint if
$a[1]\cap b[1]=\emptyset$ for every $a,b\in\mathcal{F}$, $a\neq
b$.

\begin{defn}
Let $X$ be a compact space, and $n$ a natural number. We say that
$X$ has property $I_n^\ast$ if for every $n$ many Knaster-disjoint
families $\mathcal{F}_1,\ldots,\mathcal{F}_n$ of pseudoclopens,
and every infinite disjoint family $\mathcal{F}_0$ of
pseudoclopens, there exist uncountable subfamilies
$\mathcal{G}_i\subset \mathcal{F}_i$ for $i=1,\ldots,n$, and an
infinite subfamily $\mathcal{G}_0\subset \mathcal{F}_0$ such that
for every choice $a_i\in \mathcal{G}_i$ for $i=0,\ldots,n$ we have
$a_0[0]\cap\cdots\cap a_n[0]=\emptyset$.
\end{defn}

\begin{prop}\label{imageInast}
Let $X$ be a compact space with property $I_n^\ast$, and $f:X\To
Y$ a continuous onto map. Then $Y$ has property $I_n^\ast$.
\end{prop}

Proof: Analogous to Proposition~\ref{imageIn}.$\qed$

\begin{prop}\label{onecccInast}
Let $X_0,\ldots,X_n$ be infinite compact spaces such that
$X=X_0\times\cdots\times X_n$ has poperty $I_n^\ast$. Then there
exists $i,j\leq n$, $i\neq j$ such that $X_i$ and $X_j$ satisfy
Knaster's condition.
\end{prop}

Proof: This is equivalent to say for every $i\leq n$ there exists
$j\neq i$ such that $X_j$ satisfies Knaster's condition. We prove
this statement for $i=0$. By contradiction, if this was false,
then we can find an uncountable Knaster-disjoint family
$\mathcal{F}_j$ of strongly nonempty pseudoclopens in $X_j$, for
every $j=1,\ldots,n$. Consider also $\mathcal{F}_0$ an infinite
disjoint family of pseudoclopens of $X_0$. In a similar way as we
did in Proposition~\ref{onecccIn}, these families can be lifted to
families of pseudoclopens of $X$ that violate property
$I_n^\ast$.$\qed$

\begin{defn}
A family $\mathcal{B}_0$ of pseudoclopens of $X$ will be called a
strong basis if for every pseudoclopen $a$ of $X$ there exists
$b\in\mathcal{B}_0$ such that $b\prec a$.
\end{defn}

Notice that if $X$ is a totally disconnected compact space, then
the family $\mathcal{B}_0$ of pseudoclopen sets of the form
$(c,c)$, $c$ clopen, constitutes a strong basis.

\begin{lem}\label{basislemmaast}
Let $\mathcal{B}$ and $\mathcal{B}_0$ be a basis and strong basis
of pseudoclopens of $X$ respectively. Then the following condition
is sufficient for $X$ having property $I_n^\ast$: For every $n$
many Knaster-disjoint families
$\mathcal{F}_1,\ldots,\mathcal{F}_n$ of pseudoclopens from
$\mathcal{B}$ and every infinite disjoint family $\mathcal{F}_0$
of pseudoclopens from $\mathcal{B}_0$, there exist uncountable
subfamilies $\mathcal{G}_i\subset \mathcal{F}_i$ for
$i=1,\ldots,n$ and an infinite subfamily $\mathcal{G}_0\subset
\mathcal{F}_0$ such that for every choice $a_i\in G_i$ for
$i=0,\ldots,n$ we have $a_0[0]\cap\cdots\cap a_n[0]=\emptyset$.
\end{lem}

Proof: Analogous to Lemma~\ref{basislemma}. Just note that we need
$\mathcal{B}_0$ to be a strong basis and not just a basis, because
when dealing with infinite instead of uncountable families, it is
not possible to fix the length of finite unions by passing to a
further subfamily. $\qed$

\section{Linearly ordered spaces}

\begin{lem}\label{clean}
Let $\mathcal{F}$ be a Knaster-disjoint family of sets. Then,
there exists at most countably many elements $a\in\mathcal{F}$
such that $\{b\in\mathcal{F} : a\cap b=\emptyset\}$ is countable.
\end{lem}

Proof: Suppose by contradiction that there are uncountably many
such elements. Then it is possible to construct by induction a
transfinite sequence $\{a_\alpha : \alpha<\omega_1\}\subset
\mathcal{F}$ of such elements such that $a_\alpha\cap
a_\beta\neq\emptyset$ for all $\alpha<\beta<\omega_1$. This
contradicts that $\mathcal{F}$ is Knaster-disjoint.$\qed$

\begin{thm}\label{orderedIn}
Every linearly ordered compact space $L$ has property $I_1$
\end{thm}

Proof: Every compact linearly ordered space $L$ is the continuous
image of a compact linearly ordered totally disconnected space
(one can consider the lexicographical product $L\times \{0,1\}$.
Therefore, we can suppose that $L$ is totally disconnected. By
Lemma~\ref{basislemma}, we have to show that for any
$\mathcal{F}_0$ and $\mathcal{F}_1$ Knaster-disjoint families of
clopen intervals of $L$, there are further uncountable subfamilies
for which all crossed intersections are empty. By
Lemma~\ref{clean}, we can suppose that each element of
$\mathcal{F}_i$ is disjoint from uncountably many elements of
$\mathcal{F}_i$, $i=0,1$.
 Notice that
Knaster-disjoint families are point-countable, that is, every
element of $L$ belongs to at most countably many intervals from
$\mathcal{F}_i$. Suppose that some interval $a\in \mathcal{F}_0$
intersects uncountably many intervals from $\mathcal{F}_1$. Except
those which contain some of the two extremes of $a$, which are at
most countably many, the rest are actually contained in $a$. In
this case it is enough to take $\mathcal{G}_0 = \{b\in
\mathcal{F}_0 : b\cap a=\emptyset\}$ and
$\mathcal{G}_1=\{b\in\mathcal{F}_1 : b\subset a\}$. The remaining
case is that every element of $\mathcal{F}_0$ intersects at most
countably many elements from $\mathcal{F}_1$ and vice-versa. In
this case we can produce by induction two $\omega_1$-sequences
$(a_\alpha)_{\alpha<\omega_1}\subset\mathcal{F}_0$ and
$(b_\alpha)_{\alpha<\omega_1}\subset\mathcal{F}_1$ with
$a_\alpha\cap b_\beta=\emptyset$ for all $\alpha,\beta<\omega_1$.
$\qed$

\begin{thm}\label{orderedInast} If $L_1,\ldots,L_n$ are linearly ordered
compact spaces, then $K=L_1\times\cdots\times L_n$ has property
$I_n^\ast$.
\end{thm}

Proof: Again, we assume that the spaces $L_j$ are totally
disconnected. Using Lemma~\ref{basislemmaast}, let $\mathcal{F}_0$
be a countably infinite disjoint family of clopens of $K$ and
$\mathcal{F}_1,\ldots,\mathcal{F}_n$ be Knaster-disjoint families
of clopen boxes of $K$ (by a box we mean a product of clopen
intervals).

$$\mathcal{F}_i = \{a^\alpha[i]=a^\alpha_1[i]\times\cdots\times a_n^\alpha[i] :
\alpha<\omega_1\}, i>0$$

 By Lemma~\ref{Knasterrectangles}
we can suppose that for every $i$ there exists $j(i)$ such that
$\{a^\alpha_{j(i)}[i] : \alpha<\omega_1\}$ is Knaster-disjoint. We
can actually assume that the map $i\mapsto j(i)$ is a bijection,
otherwise if there existed $i\neq i'$ with $j(i)=j(i')=j$ we would
be done by applying that $L_j$ has property $I_1$. After
relabelling we suppose that each family $\mathcal{H}_i =
\{a^\alpha_i[i] : \alpha<\omega_1\}$ is a Knaster-disjoint family
of clopen intervals of $L_i$. We consider two cases

Case 1: There exists $c\in \mathcal{F}_0$ and a box
$b_1\times\cdots\times b_n\subset c$ such that $A_i=
\{\alpha<\omega_1 : a^\alpha_i[i]\subset b_i\}$ is uncountable for
all $i=1,\ldots,n$. In this case, we can take $\mathcal{G}_0 =
\mathcal{F}_0\setminus \{c\}$ and $\mathcal{F}_i = \{a^\alpha[i] :
\alpha\in A_i\}$.

 Case 2: The previous case does not hold, so for
every $c\in \mathcal{F}_0$ and every box $b=b_1\times\cdots\times
b_n\subset c$ there exists $i=i(b)\in\{1,\ldots,n\}$ such that the
set $A_{i(b)}(b)= \{\alpha<\omega_1 : a^\alpha_i[i]\subset b_i\}$
is countable. Actually the set $A'_i(b)= \{\alpha<\omega_1 :
a^\alpha_i[i]\cap b_i\neq\emptyset\}$ is also countable, because
the family $\mathcal{H}_i = \{a^\alpha_i[i] : \alpha<\omega_1\}$
is point-countable, so only countably many intervals from it can
hit one of the two extremes of $b_i$. Every clopen
$c\in\mathcal{F}_0$ is a finite union of boxes, so there is a
countable family $\mathcal{B}$ of boxes contained in elements of
$\mathcal{F}_0$ such that every element of $\mathcal{F}_0$ is a
finite union of elements of $\mathcal{B}$. The subfamilies
$\mathcal{G}_0 = \mathcal{F}_0$ and $\mathcal{G}_i = \{a^\alpha[i]
: \alpha<\omega_1, \alpha\not\in\bigcup_{b\in\mathcal{B}}
A'_{i(b)}(b)\}$ satisfy the desired properties.$\qed$

In order to complete the proof of
Theorem~\ref{orderindecomposable}, after
Theorem~\ref{orderedInast}, Proposition~\ref{imageInast} and
Proposition~\ref{onecccInast}, it only remains to pass from
Knaster's condition to separability. It is a classical result of
Knaster \cite{Knaster} that a linearly ordered space which
satisfies Knaster's condition is separable. We need just a little
bit more:

\begin{prop}
Let $L_1,\ldots,L_n$ be linearly ordered compact spaces and
$f:L_1\times\cdots\times L_n\To X$ an onto continuous map. If $X$
satisfies Knaster's condition, then $X$ is separable.
\end{prop}

Proof: Let $Y\subset L_1\times\cdots\times L_n$ be a closed
subspace such that $f:Y\To X$ is an irreducible onto map. Recall
that a continuous function is irreducible if $f(Y')\neq X$
whenever $Y'\neq Y$ is a proper closed subset of $Y$. A standard
argument using Zorn's lemma yields the existence of such a $Y$.
Knaster's condition is preserved by irreducible preimages, so $Y$
satisfies Knaster's condition (to see this, associate to an
uncountable family $\mathcal{F}$ of nonempty open sets of $Y$, the
family $\{X\setminus f(X\setminus U) : U\in\mathcal{F}\}$ of
nonempty open sets of $X$). Now let $p_i:L_1\times\cdots\times
L_n\To L_i$ be the natural projection, and $K_i=p_i(Y)$. Then
$p_i(Y)$ is a linearly ordered compact with Knaster's property, so
by Knaster's result \cite{Knaster} it is separable. On the other
hand, $Y\subset p_1(Y)\times\cdots\times p_n(Y)$, so
$f(p_1(Y)\times\cdots\times p_n(Y))= X$, therefore $X$ is also
separable.$\qed$

\section{Spaces of finite sets}\label{Sectionsigmaens}

For a natural number $n$ and an uncountable set $\Gamma$, let
$\sigma_n(\Gamma)$ denote the family of subsets of $\Gamma$ of
cardinality less than or equal to $n$. This is a closed subset of
$2^\Gamma$, so we view $\sigma_n(\Gamma)$ as a compact topological
space. A basis for its topology are the sets of the form
$$\Phi_A^B = \{ C\in\sigma_n(\Gamma) : A\subset C\subset\Gamma\setminus B\},$$ where $A$ and $B$ are finite subsets of $\Gamma$.

The topological classification of the spaces which are finite or
countable products of spaces $\sigma_n(\Gamma)$ is studied in
\cite{AvilesCountableProducts}. In the case of finite products, in
which we are interested now, if
$\sigma_1(\Gamma)^{e_1}\times\cdots\times \sigma_n(\Gamma)^{e_n}$
is homeomorphic to $\sigma_1(\Gamma)^{f_1}\times\cdots\times
\sigma_n(\Gamma)^{f_n}$, where $n$ and each $e_i$, $f_i$ are
natural numbers, then $e_i=f_i$ for every $i$.

In this section, we will determine when a finite product of spaces $\sigma_n(\Gamma)$ can be mapped continuously onto another. An obvious sufficient condition for the existence of a continuous onto map between such finite products is the following:

\begin{lem}\label{sigmaenslemma}
Let $(n_1,\ldots,n_r)$ and $(m_1,\ldots,m_s)$ be two finite sets
of natural numbers. Suppose that there exist sets $S_i\subset
\{1,\ldots,r\}$ for $i=1,\ldots,s$, which are pairwise disjoint
and such that $m_i\leq \sum_{j\in S_i} n_j$ for every
$i\in\{1,\ldots,s\}$. Then the space
$\sigma_{n_1}(\Gamma)\times\cdots\times\sigma_{n_r}(\Gamma)$ maps
continuously onto
$\sigma_{m_1}(\Gamma)\times\cdots\times\sigma_{m_s}(\Gamma)$.
\end{lem}

Proof: The first remark is that $\sigma_n(\Gamma)$ maps
continuously onto $\sigma_m(\Gamma)$ if $m\leq n$. Namely, fix
$\gamma_0\in \Gamma$ and then define $f:\sigma_n(\Gamma)\To
\sigma_{n-1}(\Gamma\setminus\{\gamma_0\})$, by
$f(x)=x\setminus\{\gamma_0\}$ if $\gamma_0\in x$, and
$f(x)=\emptyset$ if $\gamma_0\not\in x$. The second remark is the
existence, for a finite set $S$ of natural numbers whose sum is
$\Sigma(S)$, of the union map $$u:\prod_{n\in S}
\sigma_{n}(\Gamma)\To \sigma_{\Sigma(S)}(\Gamma),$$ $u(x) =
\bigcup_{n\in S}x_n$. The two remarks together provide that
$\prod_{n\in S}\sigma_{n}(\Gamma)$ maps onto $\sigma_m(\Gamma)$
whenever $m\leq \Sigma(S)$. The proof of the lemma is obtained by
applying this fact to $S=S_i$ for every $i$, and considering
product maps.$\qed$

In Theorem \ref{sigmaensimage} below we shall prove that the
sufficient condition of the previous lemma is actually necessary.
Indeed, we shall obtain stronger indecomposibility properties. An
$m$-point family of sets is a family $\mathcal{F}$ such that every
subfamily of cardinality $m+1$ has empty intersection.

\begin{defn}
Let $m_\ast = (m_1,\ldots,m_s)$ be a finite sequence of natural
numbers. We say that a compact space has property $I[m_\ast]$ if
for every $\mathcal{F}_1\ldots,\mathcal{F}_s$ uncountable families
of clopen sets such that $\mathcal{F}_i$ is an $m_i$-point familiy
for every $i$, then there exists uncountable subfamilies
$\mathcal{G}_i\subset \mathcal{F}_i$ such that
$\bigcup_1^s\mathcal{G}_i$ is a $((\sum_1^s m_i)-1)$-point family.
\end{defn}

\begin{thm}\label{sigmaensimage}
Let $(n_1,\ldots,n_r)$ and $(m_1,\ldots,m_s)$ be two finite sets
of natural numbers. The following are equivalent:

\begin{enumerate}
\item
$K=\sigma_{n_1}(\Gamma)\times\cdots\times\sigma_{n_r}(\Gamma)$
does not have property $(m_1,\ldots,m_s)$.

\item The space
$\sigma_{n_1}(\Gamma)\times\cdots\times\sigma_{n_r}(\Gamma)$ maps
continuously onto
$\sigma_{m_1}(\Gamma)\times\cdots\times\sigma_{m_s}(\Gamma)$.

\item There exist sets $S_i\subset \{1,\ldots,r\}$ for $i=1,\ldots,s$, which are disjoint and such that $m_i\leq \sum_{j\in S_i} n_j$ for every
$i\in\{1,\ldots,s\}$.

\end{enumerate}
\end{thm}

Proof: $[3\Rightarrow 2]$ is Lemma \ref{sigmaenslemma}, $[2\Rightarrow
1]$ is clear: Consider $g:\sigma_{n_1}(\Gamma)\times\cdots\times\sigma_{n_r}(\Gamma)\To \sigma_{m_1}(\Gamma)\times\cdots\times\sigma_{m_s}(\Gamma)$ onto, and the families
$$\mathcal{F}_i = \{ \{x : \gamma\in g(x)_i\} : \gamma\in\Gamma\}, \ i=1,\ldots,s.$$ These families witness the failure of $I[m_\ast]$.
It remains to prove $[1\Rightarrow 3]$. As $K$ does not have
property $(m_1,\ldots,m_s)$, there exist families
$\mathcal{F}_1,\ldots,\mathcal{F}_s$ in $K$ such that
$\mathcal{F}_i$ is an $m_i$-point family, but for any uncountable
subfamilies $\mathcal{G}_i\subset \mathcal{F}_i$, the union
$\bigcup_1^s\mathcal{G}_i$ is not a $(\sum_1^s m_i-1)$-point
family.

For a fixed $i$, we can suppose that each clopen $x\in
\mathcal{F}_i$ is of the form
$$ (\star)\ \ x = \bigcup_{p=1}^{k(i)}\Phi_{a(x,p,1)}^{b(x,p,1)}\times\cdots\times\Phi_{a(x,p,r)}^{b(x,p,r)},$$

where $a(x,p,j)$ and $b(x,p,j)$ are finite subsets of $\Gamma$,
$\Phi_{a(x,p,j)}^{b(x,p,j)}$ is a basic clopen subset of
$\sigma_{n_j}(\Gamma)$. Moreover we can suppose that for fixed
$i,p,j$, $\{a(x,p,j) : x\in \mathcal{F}_i\}$ and $\{b(x,p,j) :
x\in \mathcal{F}_i\}$ form $\Delta$-systems of constant
cardinality with roots $A(i,p,j)$ and $B(i,p,j)$. We write
$a(x,p,j) = A(i,p,j) \cup \alpha(x,p,j)$ and  $b(x,p,j) = B(i,p,j)
\cup \beta(x,p,j)$ separating the root and the disjoint part of
the $\Delta$-system. We also write $|\alpha|(i,p,j) =
|\alpha(x,p,j)|$, $x\in \mathcal F_i$. By passing to further
uncountable subfamilies we can also assume that $$(\star\star)
(\alpha(x,p,j)\cup \beta(x,p,j)) \cap (a(x', p',j')\cup
b(x',p',j'))=\emptyset$$ whenever $x\neq x'$.

Because each $\mathcal F_i$ is an $m_i$-point family but the
family $\bigcup\mathcal F_i$ is not a $(\sum_1^s m_i-1)$-point
family, there must exist $x^i_1,\ldots,x^i_{m_i}\in \mathcal F_i$
such that $\bigcap_{i=1}^r\bigcap_{q=1}^{m_i}x^i_q \neq
\emptyset$.

Each element $x\in S_i$ is of the form $(\star)$ so it follows
that there exist also $p_i^1,\ldots, p_i^{m_i}$ for every $i$ such
that
$$\bigcap_{i=1}^r\bigcap_{q=1}^{m_i}\Phi_{a(x_i^q,p_i^q,1)}^{b(x_i^q,p_i^q,1)}\times\cdots\times\Phi_{a(x_i^q,p_i^q,r)}^{b(x_i^q,p_i^q,r)}
\neq \emptyset$$

We define the sets $S_i$ appearing in statement (3) of the Theorem in the following way:

$$S_i = \left\{ j\in\{1,\ldots,r\} : \exists \bar{q} :  |\alpha|(i,p_i^{\bar{q}},j)+\left|\bigcup_{q=1}^{m_i} a(x_i^q,p_i^q,j)\right| > n_j \right\}$$

Claim 1: $S_i\cap S_{i'} =\emptyset$.

Proof: Suppose $j\in S_i \cap S_{i'}$. Consider witnesses
$\bar{q}$ and $\bar{q}'$ for $j\in S_i$ and $j\in S_{i'}$
respectively. Assume that $|\alpha|(i,p_i^{\bar{q}},j) \leq
|\alpha|(i',p_{i'}^{\bar{q}'},j)$. We know that

$$ \left(\bigcap_{q=1}^{m_i} \Phi_{a(x_i^q,p_i^q,j)}^{b(x_i^q,p_i^q,j)}\right)\cap \left(\bigcap_{q=1}^{m_{i'}} \Phi_{a(x_{i'}^q,p_{i'}^q,j)}^{b(x_{i'}^q,p_{i'}^q,j)}\right)\neq\emptyset.  $$

In particular, $$\left|a(x_{i'}^{\bar{q}'},p_{i'}^{\bar{q}'},j) \cup \bigcup_{q=1}^{m_i}a(x_i^q,p_i^q,j)\right|\leq n_j,$$

but \begin{eqnarray*} n_j &<&
|\alpha|(i,p_i^{\bar{q}},j)+\left|\bigcup_{q=1}^{m_i}
a(x_i^q,p_i^q,j)\right| \leq
|\alpha|(i',p_{i'}^{\bar{q}'},j)+\left|\bigcup_{q=1}^{m_i}
a(x_i^q,p_i^q,j)\right|\\ &\leq&
\left|a(x_{i'}^{\bar{q}'},p_{i'}^{\bar{q}'},j) \cup
\bigcup_{q=1}^{m_i}a(x_i^q,p_i^q,j)\right|,\end{eqnarray*}  a
contradiction.

Claim 2: Fix $i\in\{1,\ldots,s\}$. For every $\bar{q}\in
\{1,\ldots,m_i\}$ there exists $j\in\{1,\ldots,r\}$ such that
$\left|\bigcup_{q=1}^{m_i} a(x,p_i^q,j)\right| +
|\alpha|(p_i^{\bar{q}},i,j) > n_j$.

Proof: Given $x\in \mathcal F_i\setminus\{x_i^q :
q=1,\ldots,m_i\}$, because $\mathcal{F}_i$ is an $m_i$-point
family,
$$
\Phi_{a(x,p_i^{\bar{q}},1)}^{b(x,p_i^{\bar{q}},1)}\times\cdots\times\Phi_{a(x,p_i^{\bar{q}},r)}^{b(x,p_i^{\bar{q}},r)}\cap
\bigcap_{q=1}^{m_i}\Phi_{a(x_i^q,p_i^q,1)}^{b(x_i^q,p_i^q,1)}\times\cdots\times\Phi_{a(x_i^q,p_i^q,r)}^{b(x_i^q,p_i^q,r)}
= \emptyset,$$ thus there exists $j$ such that
$$ \Phi_{a(x,p_i^{\bar{q}},j)}^{b(x,p_i^{\bar{q}},j)}\cap \bigcap_{q=1}^{m_i}\Phi_{a(x_i^q,p_i^q,j)}^{b(x_i^q,p_i^q,j)} = \emptyset.$$

Since we know that $\bigcap_{q=1}^{m_i}\Phi_{a(x_i^q,p_i^q,j)}^{b(x_i^q,p_i^q,j)} \neq \emptyset$, there are only three possibilities that such an intersection of basic clopen sets of the form $\Phi_F^G$ is empty in $\sigma_{n_j}(\Gamma)$: \begin{enumerate}
\item either $\left|a(x,p_i^{\bar{q}},j)\cup  \bigcup_{q=1}^{m_i} a(x_i^q,p_i^q,j)\right|>n_j$,
\item or $a(x,p_i^{\bar{q}},j)\cap  \bigcup_{q=1}^{m_i} b(x_i^q,p_i^q,j)\neq\emptyset$,
\item or $b(x,p_i^{\bar{q}},j)\cap  \bigcup_{q=1}^{m_i} a(x_i^q,p_i^q,j)\neq\emptyset$.
\end{enumerate}

The first option leads immediately to the desired conclusion. The
second and third options cannot occur. The reason is that we
assumed the $\Delta$-systems to satisfy the disjointness property
$(\star\star)$. Thus, for instance if (2) happened we would have
$A(i,p_i^{\bar{q}},j)\cap B(i,p_i^{\hat{q}},j)\neq\emptyset$ for
some $\hat{q}$, which implies that
$\Phi_{a(x_i^{\bar{q}},p_i^{\bar{q}},j)}^{b(x_i^{\bar{q}},p_i^{\bar{q}},j)}\cap
\Phi_{a(x_i^{\hat{q}},p_i^{\hat{q}},j)}^{b(x_i^{\hat{q}},p_i^{\hat{q}},j)}
= \emptyset$, which is a contradiciton.

Claim 3: $\sum_{j\in S_i}n_j\geq m_i$.

Proof: Consider the function $J:\{1,\ldots,m_i\}\To
\{1,\ldots,r\}$ which associates to every $\bar{q}$ an element
$j=J(\bar{q})$ as in claim 2. It is enough to notice that
$|J^{-1}(j)|\leq n_j$ for every $j$. Namely,
$$\bigcap_{q\in
J^{-1}(j)}\Phi_{a(x_i^q,p_i^q,j)}^{b(x_i^q,p_i^q,j)}\neq
\emptyset,$$ so $|\bigcup_{\bar{q}\in
J^{-1}(j)}a(x_i^q,p_i^q,j)|\leq n_j$. But on the other hand,
$|\alpha|(i,p_i^{\bar{q}},j)>0$ for every $\bar{q}\in J^{-1}(j)$,
so $|J^{-1}(j)|\leq n_j$.$\qed$

\section{Remarks}


\begin{rem} Theorem~\ref{Ballindecomposable} is in a sense best possible,
because $B^n$ maps continuously onto $B^n\times [0,1]^\omega$.
This is a consequence of a result of Kalenda~\cite{Kalenda} that
$B_+(\Gamma)\approx P(\sigma_1(\Gamma))$ is homeomorphic to
$B_+(\Gamma)\times [0,1]$ together with the facts that $B_+$ and
$B$ are continuous image of each other, and $[0,1]^\omega$ is a
continuous image of $[0,1]$.
\end{rem}

\begin{rem} Marde\v{s}i\'{c}'s conjecture is also best possible. For example, any
non metrizable non scattered compact space maps continuously onto
$[0,1]^\omega$. Hence we can easily get continuous maps
$L_1\times\cdots\times L_n\To L_1\times\cdots\times L_{n-1}\times
[0,1]^\omega$.
\end{rem}

\begin{rem} The reader might wonder why we deal with Knaster-disjoint
families instead of simply disjoint families. The reason is that
that Knaster-disjointness behaves better with respect to products.
For instance, if we define similar properties to $I_n$ or
$I_n^\ast$ with disjoint families instead of Knaster-disjoint,
then the proofs of Proposition~\ref{productivity} and
Theorem~\ref{orderedInast} do not work any more, unless we assume
that certain colorings of the uncountable have uncountable
monochromatic sets: In the first case this can be overcome if the
compact spaces $X$ and $Y$ satisfy property (Q) of
Bell~\cite{BellRamsey}, and in the second case one may need to
assume that Suslin lines do not exist.
\end{rem}

\begin{rem}
Given an uncountable regular cardinal $\aleph$, we can define
indecomposability properties $I_n(\aleph)$, $I_n^\ast(\aleph)$ or
$I[m_\ast](\aleph)$ in a similar way but substituting ``\emph{for
any uncountable families... there exist uncountable subfamilies}''
by ``\emph{for any families of cardinality $\aleph$... there exist
subfamilies of cardinality $\aleph$}''. All the results in this
note can be rewritten in this more general way, and in particular
the Theorems \ref{Ballirreducible}, \ref{orderedIn},
\ref{orderedInast} and \ref{sigmaensimage} hold for these
properties relative to $\aleph$.
\end{rem}

\begin{rem}
After Rudin's result \cite{Rudin} that every monotonically normal
compact space is the continuous image of a linearly ordered space,
in Theorem~\ref{orderindecomposable} the assumption that the
spaces $L_i$ are linearly ordered can be substituted by the
assumption that the spaces $L_i$ are monotonically normal.
\end{rem}

\end{document}